\documentclass[graybox]{svmult}
\usepackage{mathptmx}       
\usepackage{helvet}         
\usepackage{courier}        
\usepackage{type1cm}        
\usepackage{makeidx}         
\usepackage{graphicx}        
\usepackage{multicol}        
\usepackage[bottom]{footmisc}

\usepackage{amssymb}
\usepackage{amsmath}

\begin{document}

\title{The discrete octonionic Stokes' formula revisited}
\author{Rolf S\"oren Krau\ss har, Anastasiia Legatiukand and Dmitrii Legatiuk \\[4mm] \emph{This paper is devoted to our friend and colleague Paula Cerejeiras in occasion of her 60th birthday.}}
\authorrunning{Rolf S\"oren Krau\ss har, Anastasiia Legatiuk and Dmitrii Legatiuk}
\institute{Rolf S\"oren Krau\ss har \at Universit\"at Erfurt, Germany \and Anastasiia Legatiuk \at Universit\"at Erfurt, Germany \and Dmitrii Legatiuk \at Universit\"at Erfurt, Germany  \email{dmitrii.legatiuk@uni-erfurt.de}}  
%
%
\maketitle

\abstract{In a previous work we made an attempt to set up a discrete octonionic Stokes' formula. Due to an algebraic property that we have not considered in that attempt, the formula however turned out to involve an associator term in addition to a change of sign that we already observed earlier. This associator term has an impact on the final result. In this paper we carefully revise this discrete Stokes formula taking into account this additional term. In fact the result that we now obtain is much more in line with the results that one has in the continuous setting.}


\section{Introduction}

The octonions received an enormously increasing interest from physicists and mathematicians during the last years. This is party due to the fact that the octonionic configuration seems to offer a more suitable model for extending the standard model of particle physics in the way offering also a possibility to incorporate some gravitational effects. This is a strong reason to continue also the development of analytic, in particular, function theoretical tools in the octonionic setting more intensively. However, the use of function-theoretic tools in practical applications involving PDEs typically requires a discretisation of a problem. Hence, a discrete version of differential operators needs to be considered as well. This approach gives rise to the discrete counterpart of the classical continuous theory, namely the {\itshape discrete octonionic analysis}, where the discrete version of the continuous Cauchy-Riemann operator plays the central role. During the last two years we basically considered three different kinds of discrete Cauchy-Riemann operators. In \cite{Krausshar_1}, we discretised the octonionic Cauchy-Riemann operator in terms of discrete forward and backward operators. After that, the central difference operator has been used in \cite{Krausshar_2}. The Weyl calculus perspective on the discrete Cauchy-Riemann operator has been studied in \cite{Krausshar_3,Krausshar_5} for unbounded and bounded domains. The Weyl calculus setting has the advantage that the classical star Laplacian can be factorised by a pair of Cauchy-Riemann operators. The discretisation based on backward, forward, or central difference operators do not have this property.\par 
In \cite{Krausshar_1}, the main intention was to develop a Stokes formula for the octonionic discrete setting using forward and backward operators. Note that a Stokes formula serves as a fundamental starting point to subsequentially develop generalisations of the Borel-Pompeiu formula and the Cauchy formula which in turn allows one to study Hardy spaces in the discrete setting. We already observed in \cite{Krausshar_1} correctly that the Stokes formula involve different minus signs as in the quaternionic or Clifford analytic setting. This change in sign is caused by the non-associativity of the octonionic multiplication, because we have $(\mathbf{e}_i \mathbf{e}_j)\mathbf{e}_k = -\mathbf{e}_i (\mathbf{e}_j \mathbf{e}_k)$, for mutually distinct and non-zero $i,j,k$ and if furthermore $\mathbf{e}_{i}\mathbf{e}_{j}\neq \pm \mathbf{e}_k$. Contrarily, if $i,j,k$ are not mutually distinct, or if one of these units is the real unit $1$, as well as if $\mathbf{e}_{i}\mathbf{e}_{j}= \pm \mathbf{e}_k$, then due to Artin's theorem we have instead $(\mathbf{e}_i \mathbf{e}_j)\mathbf{e}_k = +\mathbf{e}_i (\mathbf{e}_j \mathbf{e}_k)$. This aspect has not been taken into account in our original paper \cite{Krausshar_1}. Therefore, the Stokes formulae presented in \cite{Krausshar_1,Krausshar_2,Krausshar_3} need to be adapted incorporating this feature. The fact that some of these combinations lead to a plus sign implies that an associator term appears in the discrete Stokes' formula and consequently also in the Borel-Pompeiu formula. In fact, this produces a closer analogy to the continuous case where also an associator appears. We hereby fill in an important gap in our preceding research, since we were now able to understand better how the discrete setting should converge to the continuous case where the associator appears.\par 
This short paper is entirely devoted to the setting of forward and backward difference operators. A similar correction needs to be applied to the approach using the central difference operator, which will be published in our follow-up paper. This forthcoming paper also discusses the effect of the associator term on the results regarding the Hardy spaces discussed in \cite{Krausshar_2}.\par 


\section{Preliminaries}

The basic context is the $8$-dimensional Euclidean vector space $\mathbb{R}^{8}$, where the standard vectors are denoted by $\mathbf{e}_{k}$, $k=0,1,\ldots,7$. A vector from $\mathbb{R}^8$ can be expressed as usual in terms of its real coordinates in the way  $\mathbf{x}=(x_{0}, x_{1},\ldots, x_{7})$. Next, $\mathbb{R}^8$-vectors can also be described as octonions 
\begin{equation*}
x = x_{0}\mathbf{e}_{0}+x_{1}\mathbf{e}_{1}+x_{2}\mathbf{e}_{2}+x_{3}\mathbf{e}_{3}+x_{4}\mathbf{e}_{4}+x_{5}\mathbf{e}_{5}+x_{6}\mathbf{e}_{6}+x_{7}\mathbf{e}_{7},
\end{equation*}
where we now additionally identify $\mathbf{e}_{4}=\mathbf{e}_{1}\mathbf{e}_{2}$, $\mathbf{e}_{5}=\mathbf{e}_{1}\mathbf{e}_{3}$, $\mathbf{e}_{6}=\mathbf{e}_{2}\mathbf{e}_{3}$ and \linebreak $\mathbf{e}_{7}=\mathbf{e}_{4}\mathbf{e}_{3}=(\mathbf{e}_{1}\mathbf{e}_{2})\mathbf{e}_{3}$. Moreover, we have $\mathbf{e}_{i}^{2}=-1$ and  $\mathbf{e}_{0}\mathbf{e}_{i}=\mathbf{e}_{i}\mathbf{e}_{0}$ for all $i=1,\ldots,7$, and $\mathbf{e}_{i}\mathbf{e}_{j}=-\mathbf{e}_{j}\mathbf{e}_{i}$ for all mutual distinct $i,j\in\left\{1,\ldots,7\right\}$, as well as $\mathbf{e}_{0}$ is the neutral element and, therefore, often will be omitted. This definition endows $\mathbb{R}^8$ additionally with a multiplicative closed structure. The octonionic multiplication is closed but not associative, for instance we have $(\mathbf{e}_{i}\mathbf{e}_{j})\mathbf{e}_{k}=-\mathbf{e}_{i}(\mathbf{e}_{j}\mathbf{e}_{k})$ for mutually distinct and non-zero indices $i,j,k$ satisfying additionally $\mathbf{e}_{i}\mathbf{e}_{j}\neq \pm \mathbf{e}_k$. Otherwise, the multiplication is associative $(\mathbf{e}_{i}\mathbf{e}_{j})\mathbf{e}_{k}=\mathbf{e}_{i}(\mathbf{e}_{j}\mathbf{e}_{k})$.\par
Let us consider the unbounded uniform lattice $h\mathbb{Z}^{8}$ with the lattice constant $h>0$, which is defined in the classical way as follows
\begin{equation*}
h \mathbb{Z}^{8} :=\left\{\mathbf{x} \in {\mathbb{R}}^{8}\,|\, \mathbf{x} = (m_{0}h, m_{1}h,\ldots, m_{7}h), m_{j} \in \mathbb{Z}, j=0,1,\ldots,7\right\}.
\end{equation*}
Next, we define the classical forward and backward differences $\partial_{h}^{\pm j}$ as
\begin{equation}
\label{Finite_differences}
\begin{array}{lcl}
\partial_{h}^{+j}f(mh) & := & h^{-1}(f(mh+\mathbf{e}_jh)-f(mh)), \\
\partial_{h}^{-j}f(mh) & := & h^{-1}(f(mh)-f(mh-\mathbf{e}_jh)),
\end{array}
\end{equation}
for discrete functions $f(mh)$ with $mh\in h\mathbb{Z}^{8}$. In the sequel, we consider functions defined on $\Omega_{h} \subset  h\mathbb{Z}^{8}$ and taking values in octonions $\mathbb{O}$.\par
By using the finite difference operators~(\ref{Finite_differences}), we introduce a {\itshape discrete forward Cauchy-Riemann operator} $D^{+}\colon l^{p}(\Omega_{h},\mathbb{O})\to l^{p}(\Omega_{h},\mathbb{O})$ and a {\itshape discrete backward Cauchy-Riemann operators} $D^{-}\colon l^{p}(\Omega_{h},\mathbb{O})\to l^{p}(\Omega_{h},\mathbb{O})$ as follows
\begin{equation}
\label{Cauchy_Riemann_operators_discrete}
D^{+}_{h}:=\sum_{j=0}^{7} \mathbf{e}_j\partial_{h}^{+j}, \quad D^{-}_{h}:=\sum_{j=0}^{7} \mathbf{e}_j\partial_{h}^{-j}.
\end{equation}


\section{Discrete octonionic Stokes' formula revisited}

We present the following theorem:
\begin{theorem}
The discrete Stokes' formula for the whole lattice $h\mathbb{Z}^{8}$ is given by
\begin{equation}
\label{Discrete_Stokes_whole_space}
\begin{array}{c}
\displaystyle \sum_{m\in \mathbb{Z}^{8}}  \left[ \left( g(mh)D_h^{+}\right) f(mh) + g(mh) \left( D_h^{-}f(mh) \right) \right] h^8 =  \\[2mm]
\displaystyle = 2\sum_{m\in \mathbb{Z}^{8}} \sum\limits_{s=1}^{7}  \sum_{i\in I_{s}} \sum_{\stackrel{j\in I_{s}}{j\neq i}}^{7} \sum_{\stackrel{k=1}{k\notin I_{s}}}^{7} \left[g_{i}(mh)\mathbf{e}_{i}\left(\partial_{h}^{-j}\mathbf{e}_{j}f_{k}(mh)\mathbf{e}_{k}\right)\right]h^{8}
\end{array}  
\end{equation}
for all discrete functions $f$ and $g$ such that the series converge, where the index sets $I_{s}$, $s=1,\ldots,7$ are given by
\begin{equation*}
\begin{array}{cclcclcclccl}
I_{1} & := & \left\{1,2,4\right\}, & I_{2} & := & \left\{1,3,5\right\}, & I_{3} & := & \left\{1,6,7\right\}, & I_{4} & := & \left\{2,3,6\right\}, \\[2mm]
I_{5} & := & \left\{2,5,7\right\}, & I_{6} & := & \left\{3,4,7\right\}, & I_{7} & := & \left\{4,5,6\right\}.
\end{array}
\end{equation*}
\end{theorem}
\begin{proof}
The proof starts by working with the first term on the left-hand side in~(\ref{Discrete_Stokes_whole_space}) and by using the definition of the discrete forward Cauchy-Riemann operator~(\ref{Cauchy_Riemann_operators_discrete}):
\begin{equation*}
\begin{array}{rl}
& \displaystyle \sum_{m\in \mathbb{Z}^{8}}\left[ g(mh)D_h^{+}\right] f(mh) = \sum\limits_{m\in \mathbb{Z}^{8}} \sum\limits_{j=0}^{7} \left[\partial_{h}^{+j}g(mh)\mathbf{e}_{j}\right] f(mh) h^{8} \\
\\
= & \displaystyle \sum\limits_{m\in \mathbb{Z}^{8}} \sum\limits_{j=0}^{7}\sum \limits_{i=0}^{7}\sum \limits_{k=0}^{7} \left[\partial_{h}^{+j}g_{i}(mh)\mathbf{e}_{i}\mathbf{e}_{j}\right] f_{k}(mh)\mathbf{e}_{k} h^{8}.
\end{array}
\end{equation*}
The multiplication of the last expression leads to 512 summands containing products of basis elements of the form $(\mathbf{e}_{i}\mathbf{e}_{j})\mathbf{e}_{k}$. It is important to recall, that octonions $\mathbb{O}$ are not associative in $\mathbb{R}^{8}$, but always contain an associative sub-algebra. This fact results in the following multiplication rules
\begin{equation*}
\begin{array}{rcll}
(\mathbf{e}_{i}\mathbf{e}_{j})\mathbf{e}_{k} & = & -\mathbf{e}_{i}(\mathbf{e}_{j}\mathbf{e}_{k}) & \mbox{for mutually distinct, non-zero } i,j,k \mbox{ and } \mathbf{e}_{i}\mathbf{e}_{j} \neq \pm \mathbf{e}_{k},  \\[2mm]
(\mathbf{e}_{i}\mathbf{e}_{j})\mathbf{e}_{k} & = & \mathbf{e}_{i}(\mathbf{e}_{j}\mathbf{e}_{k}) & \mbox{otherwise}.
\end{array}
\end{equation*}
The associative part has been overseen in \cite{Krausshar_1,Krausshar_2,Krausshar_3}. Evidently, presenting all 512 summands is not possible here. Therefore, let us generally explain the main idea of the next intermediate steps. The non-associative part of the product contains exactly 168 terms. In order to proceed with the proof, we add and subtract these 168 terms. In that case, the associative part contains now 512 terms all written with the positive sign. Hence, these associative part can be written again by the help of a triple summation over indices $i$, $j$, and $k$ as follows
\begin{equation*}
\sum\limits_{m\in \mathbb{Z}^{8}} \sum\limits_{j=0}^{7}\sum \limits_{i=0}^{7}\sum \limits_{k=0}^{7} \partial_{h}^{+j}g_{i}(mh)\mathbf{e}_{i}\left[\mathbf{e}_{j} f_{k}(mh)\mathbf{e}_{k}\right] h^{8}.
\end{equation*}
The remaining 168 non-associative terms are grouped at first into sums of the form
\begin{equation*}
\sum_{\stackrel{k=1}{k\neq 1,2,4}}^{7} \partial_{h}^{+2}g_{1}(mh)\mathbf{e}_{1}\left(\mathbf{e}_{3}f_{k}(mh)\mathbf{e}_{k}\right)
\end{equation*}
written for all different sets of indices $i,j,k\neq 0$. This grouping leads to 42 sums, where each summation is written for $k=1,\ldots,7$ but with three terms always omitted. Further we introduce the following index sets:
\begin{equation*}
\begin{array}{cclcclcclccl}
I_{1} & := & \left\{1,2,4\right\}, & I_{2} & := & \left\{1,3,5\right\}, & I_{3} & := & \left\{1,6,7\right\}, & I_{4} & := & \left\{2,3,6\right\}, \\[2mm]
I_{5} & := & \left\{2,5,7\right\}, & I_{6} & := & \left\{3,4,7\right\}, & I_{7} & := & \left\{4,5,6\right\}.
\end{array}
\end{equation*}
These index sets allow us to write 42 sums in a short form as follows
\begin{equation*}
-2\sum_{m\in \mathbb{Z}^{8}} \sum\limits_{s=1}^{7}  \sum_{i\in I_{s}} \sum_{\stackrel{j\in I_{s}}{j\neq i}}^{7} \sum_{\stackrel{k=1}{k\notin I_{s}}}^{7} \left[g_{i}(mh)\mathbf{e}_{i}\left(\partial_{h}^{+j}\mathbf{e}_{j}f_{k}(mh)\mathbf{e}_{k}\right)\right]h^{8},
\end{equation*}
where the factor 2 comes from the fact, that we added and subtracted the non-associative terms. Thus, after all these manipulations we arrived at the following expression:
\begin{equation*}
\begin{array}{c}
\displaystyle \sum\limits_{m\in \mathbb{Z}^{8}} \sum\limits_{j=0}^{7}\sum \limits_{i=0}^{7}\sum \limits_{k=0}^{7} \partial_{h}^{+j}g_{i}(mh)\mathbf{e}_{i}\left[\mathbf{e}_{j} f_{k}(mh)\mathbf{e}_{k}\right] h^{8} - \\ \\
\displaystyle -2\sum_{m\in \mathbb{Z}^{8}} \sum\limits_{s=1}^{7}  \sum_{i\in I_{s}} \sum_{\stackrel{j\in I_{s}}{j\neq i}}^{7} \sum_{\stackrel{k=1}{k\notin I_{s}}}^{7} \left[g_{i}(mh)\mathbf{e}_{i}\left(\partial_{h}^{+j}\mathbf{e}_{j}f_{k}(mh)\mathbf{e}_{k}\right)\right]h^{8},
\end{array}
\end{equation*}
which needs to be further modified. Now we use the definition of finite difference operators~(\ref{Finite_differences}) and obtain
\begin{equation*}
\begin{array}{c}
\displaystyle \sum\limits_{m\in \mathbb{Z}^{8}} \left[ \sum\limits_{j=0}^{7}\sum \limits_{i=0}^{7}\sum \limits_{k=0}^{7} \left(g_{i}(mh+\mathbf{e}_{j}h)\mathbf{e}_{i}f_{k}(mh)\left(\mathbf{e}_{j}\mathbf{e}_{k}\right) - g_{i}(mh)\mathbf{e}_{i}f_{k}(mh)\left(\mathbf{e}_{j}\mathbf{e}_{k}\right) \right)\right.- \\ \\
\displaystyle \left. -2\sum\limits_{s=1}^{7}  \sum_{i\in I_{s}} \sum_{\stackrel{j\in I_{s}}{j\neq i}}^{7} \sum_{\stackrel{k=1}{k\notin I_{s}}}^{7} \left(g_{i}(mh+\mathbf{e}_{j}h)\mathbf{e}_{i}f_{k}(mh)\left(\mathbf{e}_{j}\mathbf{e}_{k}\right) - g_{i}(mh)\mathbf{e}_{i}f_{k}(mh)\left(\mathbf{e}_{j}\mathbf{e}_{k}\right) \right)\right]h^{8}.
\end{array}
\end{equation*}
Next, performing the change of variables in the last expression, we obtain
\begin{equation*}
\begin{array}{c}
\displaystyle \sum\limits_{m\in \mathbb{Z}^{8}} \left[ \sum\limits_{j=0}^{7}\sum \limits_{i=0}^{7}\sum \limits_{k=0}^{7} \left(g_{i}(mh)\mathbf{e}_{i}f_{k}(mh-\mathbf{e}_{j}h)\left(\mathbf{e}_{j}\mathbf{e}_{k}\right) - g_{i}(mh)\mathbf{e}_{i}f_{k}(mh)\left(\mathbf{e}_{j}\mathbf{e}_{k}\right) \right)\right.- \\ \\
\displaystyle \left. -2\sum\limits_{s=1}^{7}  \sum_{i\in I_{s}} \sum_{\stackrel{j\in I_{s}}{j\neq i}}^{7} \sum_{\stackrel{k=1}{k\notin I_{s}}}^{7} \left(g_{i}(mh)\mathbf{e}_{i}f_{k}(mh-\mathbf{e}_{j}h)\left(\mathbf{e}_{j}\mathbf{e}_{k}\right) - g_{i}(mh)\mathbf{e}_{i}f_{k}(mh)\left(\mathbf{e}_{j}\mathbf{e}_{k}\right) \right)\right]h^{8} = \\
\\
\displaystyle = \sum\limits_{m\in \mathbb{Z}^{8}} \left[- \sum\limits_{j=0}^{7}\sum \limits_{i=0}^{7}\sum \limits_{k=0}^{7} \left(g_{i}(mh)\mathbf{e}_{i}f_{k}(mh)\left(\mathbf{e}_{j}\mathbf{e}_{k}\right)-g_{i}(mh)\mathbf{e}_{i}f_{k}(mh-\mathbf{e}_{j}h)\left(\mathbf{e}_{j}\mathbf{e}_{k}\right) \right)\right.- \\ \\
\displaystyle \left. +2\sum\limits_{s=1}^{7}  \sum_{i\in I_{s}} \sum_{\stackrel{j\in I_{s}}{j\neq i}}^{7} \sum_{\stackrel{k=1}{k\notin I_{s}}}^{7} \left( g_{i}(mh)\mathbf{e}_{i}f_{k}(mh)\left(\mathbf{e}_{j}\mathbf{e}_{k}\right) - g_{i}(mh)\mathbf{e}_{i}f_{k}(mh-\mathbf{e}_{j}h)\left(\mathbf{e}_{j}\mathbf{e}_{k}\right)\right)\right]h^{8},
\end{array}
\end{equation*}
and by using the definition of the difference operator $\partial_{h}^{-j}$ this expression can be rewritten as follows
\begin{equation*}
\begin{array}{c}
\displaystyle -\sum\limits_{m\in \mathbb{Z}^{8}} \sum\limits_{j=0}^{7}\sum \limits_{i=0}^{7}\sum \limits_{k=0}^{7} g_{i}(mh)\mathbf{e}_{i}\left(\partial_{h}^{-j}\mathbf{e}_{j}f_{k}(mh)\mathbf{e}_{k}\right)h^{8} + \\ \\
\displaystyle  +2\sum\limits_{m\in \mathbb{Z}^{8}}\sum\limits_{s=1}^{7}  \sum_{i\in I_{s}} \sum_{\stackrel{j\in I_{s}}{j\neq i}}^{7} \sum_{\stackrel{k=1}{k\notin I_{s}}}^{7} g_{i}(mh)\mathbf{e}_{i}\left(\partial_{h}^{-j}\mathbf{e}_{j}f_{k}(mh)\mathbf{e}_{k}\right)h^{8}.
\end{array}
\end{equation*}
Finally, we obtain the following expression written by the help of the discrete backward Cauchy-Riemann operator:
\begin{equation*}
\begin{array}{c}
\displaystyle -\sum\limits_{m\in \mathbb{Z}^{8}}  g(mh)\left(D_{h}^{-}f(mh)\right)h^{8} + 2\sum\limits_{m\in \mathbb{Z}^{8}}\sum\limits_{s=1}^{7}  \sum_{i\in I_{s}} \sum_{\stackrel{j\in I_{s}}{j\neq i}}^{7} \sum_{\stackrel{k=1}{k\notin I_{s}}}^{7} g_{i}(mh)\mathbf{e}_{i}\left(\partial_{h}^{-j}\mathbf{e}_{j}f_{k}(mh)\mathbf{e}_{k}\right)h^{8}.
\end{array}
\end{equation*}
Bringing now the summation with $D_{h}^{-}$ to the left-hand side, we obtain~(\ref{Discrete_Stokes_whole_space}). This finishes the proof of the discrete octonionic Stokes' formula.
\end{proof}





\end{document}